\documentclass[12pt]{l4dc2021}

\usepackage{mathrsfs}
\usepackage{mathtools}
\usepackage{epstopdf}

\title[Learning Approximate Forward Reachable Sets Using Separating Kernels]{Learning Approximate Forward Reachable Sets \\ Using Separating Kernels}
\usepackage{times}

\author{%
 \Name{Adam J. Thorpe} \Email{ajthor@unm.edu} \\
 \Name{Kendric R. Ortiz} \Email{kendric@unm.edu} \\
 \Name{Meeko M. K. Oishi} \Email{oishi@unm.edu} \\
 \addr Department of Electrical and Computer Engineering, University of New Mexico%
}

\begin{document}

\maketitle

\begin{abstract}%
We present a data-driven method for computing approximate forward reachable sets using separating kernels in a reproducing kernel Hilbert space. We frame the problem as a support estimation problem, and learn a classifier of the support as an element in a reproducing kernel Hilbert space using a data-driven approach. Kernel methods provide a computationally efficient representation for the classifier that is the solution to a regularized least squares problem. The solution converges almost surely as the sample size increases, and admits known finite sample bounds. This approach is applicable to stochastic systems with arbitrary disturbances and neural network verification problems by treating the network as a dynamical system, or by considering neural network controllers as part of a closed-loop system. We present our technique on several examples, including a spacecraft rendezvous and docking problem, and two nonlinear system benchmarks with neural network controllers.%
\end{abstract}



\section{Introduction}


Reachability analysis is important in many dynamical systems for its ability to provide assurances of desired behavior, such as reaching a desired target set or anticipating potential intersection with a set known to be unsafe.
Such analysis has shown utility in a variety of safety-critical applications, including autonomous cars, UAVs, and spacecraft. 
The forward {\em stochastic reachable set} describes the set of states that the system will reach with a non-zero likelihood.
While most approaches for stochastic reachability are model-based, data-driven approaches are important when mathematical models may not exist
or may be too complex for existing numerical approaches.
In particular, the growing presence of learning-enabled components in dynamical systems warrants the development of new tools for stochastic reachability that can accommodate neural networks, look-up tables, and other model-resistant elements.
In this paper, we propose a technique for data-driven stochastic reachability analysis that provides a convergent approximation to the true stochastic reachable set.


Significant advances have been made in verification of dynamical systems with neural network components.  Recent work in \citet{seidman2020robust, weinan2017proposal} has shown that neural networks can be modeled as a nonlinear dynamical system, making them amenable in some cases to model-based reachability analysis.
Typically, these approaches presume that the neural network exhibits a particular structure, such as having a particular activation function, as in \citet{sidrane2019overt, tran2020nnv}, and exploit existing tools for forward reachability, such as \citet{althoff2015introduction, chen2013flow}.
Other approaches employ a mixed integer linear programming approach, as in \citet{dutta2017output, dutta2019sherlock, lomuscio2017approach}.
These techniques exploit Lipschitz constants or perform set propagation, however, the latter approach
becomes intractable for large-scale systems due to vertex facet enumeration.
Further, in practice, knowledge of the network structure or dynamics may not be available, or may be too complex to use with standard reachability methods.
Thus, additional tools are needed to efficiently compute stochastic reachable sets.

\begin{sloppypar}
Data-driven reachability methods provide convergent approximations of reachable sets with high confidence, but are typically unable to provide assured over-approximations.
Methods have been developed that use convex hulls in \cite{lew2020sampling} and scenario optimization in \cite{devonport2020estimating}. However, these approaches rely upon convexity assumptions which can be limiting in certain cases.
Other data-driven approaches leverage a class of machine learning techniques known as kernel methods, including Gaussian processes in \citet{devonport2020data} and support vector machines in  \citet{allen2014machine, rasmussen2017approximation}.
However, the approach in \citet{rasmussen2017approximation} can suffer from stability issues and does not provide probabilistic guarantees of convergence, and the approach in \cite{allen2014machine} requires that we repeatedly solve a nonlinear program offline to generate data for the SVM classifier.
The approach in \citet{devonport2020data}
requires use of a Gaussian process prior.
In contrast to these approaches, we propose a method that can accommodate nonlinear dynamical systems with arbitrary stochastic disturbances.
\end{sloppypar}

Our approach employs a class of kernel methods known as separating kernels to form a reachable set classifier.
By learning a set classifier in Hilbert space, we convert the problem of learning a set boundary into the problem of learning a classifier in a high-dimensional space of functions.
%
Our approach extends the work in \cite{de2014learning} to the problem of learning reachable sets for stochastic systems that overcomes the stability and convergence issues faced by existing kernel based approaches and allows for arbitrary stochastic disturbances.
Our main contribution is an application of the techniques presented in \cite{de2014learning} to the problem of learning approximate forward reachable sets and neural network verification.
Similar to other data-driven approaches, the approximation proposed here does not provide guarantees in the form of assured over- or under-approximations of the forward reachable set. However, although empirically derived, the approximation can be shown to converge in probability almost surely.

The paper is structured as follows:
Section \ref{section: problem formulation} formulates the problem.
We describe the application of kernel methods to compute approximate forward reachable sets in Section \ref{section: main results}, and discuss their convergence properties and finite sample bounds.
In Section \ref{section: numerical experiments}, we demonstrate our approach on a realistic satellite rendezvous and docking problem, as well as two neural network verification benchmarks from \citet{dutta2019sherlock}.
Concluding remarks are presented in Section \ref{section: conclusion}.



\section{Problem Formulation}
\label{section: problem formulation}

We use the following notation throughout:
Let $E$ be an arbitrary nonempty space, and denote the $\sigma$-algebra on $E$ by $\mathcal{E}$.
Let $\wp(E)$ denote the set of all subsets (the power set) of $E$.
If $E$ is a topological space \citep{ccinlar2011probability},
the $\sigma$-algebra generated by the set of all open subsets of $E$
is called the Borel $\sigma$-algebra, denoted by $\mathscr{B}(E)$.
If $E$ is countable, with $\sigma$-algebra $\wp(E)$, then it is called \emph{discrete}.
Let $(\Omega, \mathcal{F}, \mathbb{P})$ denote a probability space,
where $\mathcal{F}$ is the $\sigma$-algebra on $\Omega$ and
$\mathbb{P} : \mathcal{F} \rightarrow [0, 1]$ is a \emph{probability measure}
on the measurable space $(\Omega, \mathcal{F})$.
A measurable function $X : \Omega \rightarrow E$ is called a \emph{random variable} taking values in $(E, \mathcal{E})$. The image of $\mathbb{P}$ under $X$, $\mathbb{P}(X^{-1}A)$, $A \in \mathcal{E}$ is called the \emph{distribution} of $X$.
Let $T$ be an arbitrary set, and for each $t \in T$, let $X_{t}$ be a random variable.
The collection of random variables $\lbrace X_{t} : t \in T \rbrace$
on $(\Omega, \mathcal{F})$ is a \emph{stochastic process}.


\subsection{System Model}

Consider a general discrete-time system with dynamics given by:
\begin{equation}
    \label{eqn: system dynamics}
    x_{k+1} = f_{k}(x_{k}, u_{k}, \theta_{k}, w_{k})
\end{equation}
with state $x_{k} \in \mathcal{X} \subset \mathbb{R}^{n}$, input $u_{k} \in \mathcal{U} \subset \mathbb{R}^{m}$, parameters $\theta_{k} \in \Theta \subset \mathbb{R}^{p}$, and stochastic process $w$, comprised of the random variables $w_{k}$, $k \in [0, N]$, defined on the measurable space $(\mathcal{W}, \mathscr{B}(\mathcal{W}))$, $\mathcal{W} \subset \mathbb{R}^{q}$. 
The system evolves over a finite time horizon $k \in [0, N]$ from initial condition $x_{0} \in \mathcal{X}$, which may be chosen according to initial probability measure $\mathbb{P}_{0}$.
%
We assume that the dynamics and the structure of the disturbance are unknown, but that a finite sample of $M \in \mathbb{N}$ observations, $\lbrace x_{N} \rbrace_{i=1}^{M}$, taken i.i.d. from $\mathbb{P}_{N}$, are available.

A special case of this formulation is a feedforward neural network \citep{seidman2020robust, weinan2017proposal},
where $k$ denotes the layer of the network, $\theta_{k}$ are the network parameters, and the dimensionality of
the state, input, parameter, and disturbance spaces
depends on $k$ and can change at each layer.
For example, in each layer, $x_{k} \in \mathcal{X}_{k} \subset \mathbb{R}^{n(k)}$, where $n(k)$ is a map to $\mathbb{N}$ and is the dimensionality of the state at each layer.


\subsection{Forward Reachable Set}

Consider a discrete time stochastic system \eqref{eqn: system dynamics}.
We consider an initial condition $x_{0} \in \mathcal{X}$ and evolve \eqref{eqn: system dynamics} to compute the forward reachable set, $\mathscr{F}(x_{0})$.
%
Formally, we define the forward reachable set $\mathscr{F}(x_{0})$ as in \citet{lew2020sampling}, which is a modification of the definition in \citet{mitchell2007comparing} or \citet{rosolia2019sample} that is amenable to neural network verification problems.

\begin{definition}[{\citealp{lew2020sampling}}]
    Given an initial condition $x_{0} \in \mathcal{X}$,
    the forward reachable set $\mathscr{F}(x_{0})$ is defined as the set of all possible future states $x_{N}$ at time $N$ for a system following a given control sequence
    $\lbrace u_{0}, u_{1}, \ldots, u_{N-1} \rbrace$.
    \begin{equation}
        \label{eqn: forward reachable set}
        \mathscr{F}(x_{0}) := \lbrace x_{N} =
        f_{N-1} \circ \cdots \circ f_{0}(x_{0}, u_{0}, \theta_{0}, w_{0})
        \,|\, x_{0} \in \mathcal{X}, u_{k} \in \mathcal{U}, \theta_{k} \in \Theta, w_{k} \in \mathcal{W} \rbrace
    \end{equation}
\end{definition}
%
%
Note that $x_{N}$ is a random variable on the probability space $(\mathcal{X}, \mathscr{B}(\mathcal{X}), \mathbb{P}_{N})$,
and that the forward reachable set can also be viewed as the support of $x_{N}$ \citep{vinod2017forward},
which is the smallest closed set $\mathscr{F}(x_{0}) \subset \mathcal{X}$ such that $\mathbb{P}_{N}(x_{N} \in \mathscr{F}(x_{0})) = 1$.

We seek to compute an approximation of the forward reachable set in \eqref{eqn: forward reachable set} for a system \eqref{eqn: system dynamics} with unknown dynamics, by using the theory of reproducing kernel Hilbert spaces.
We formulate the problem as learning a classifier function $F$ in Hilbert space which describes the geometry of the forward reachable set boundary.
We presume that the forward reachable set \eqref{eqn: forward reachable set} is implicitly defined by the classifier $F$.
\begin{equation}
    \label{eqn: forward reachable set in RKHS}
    \mathscr{F}(x_{0}) = \lbrace x \in \mathcal{X} \,|\, F(x) = 1 \rbrace
\end{equation}
Hence, we seek to compute an \emph{empirical estimate} $\tilde{F}$ of $F$ in \eqref{eqn: forward reachable set in RKHS} using observations taken from the system evolution $\lbrace x_{N} \rbrace_{i=1}^{M}$.

%
Forward reachability often focuses on forming an over-approximation
that is a superset of $\mathscr{F}(x_{0})$.
However, because our proposed approach cannot provide such a guarantee, we aim to provide an approximation of the forward reachable set $\tilde{\mathscr{F}}(x_{0})$ by computing $\tilde{F}$ that is convergent in probability almost surely.
These probabilistic guarantees are important to ensure the consistency of our result and that our estimated classifier (and therefore the approximate reachable set) is close to the true result with high probability.
The main difficulty in this problem
arises from the fact that we do not have explicit knowledge of the dynamics \eqref{eqn: system dynamics} or place any prior assumptions on the structure of the uncertainty.
This makes the proposed method a useful technique in the context of existing stochastic reachability toolsets, since it can be used on black-box systems with arbitrary disturbances.
Note that because we seek to estimate a classifier $\tilde{F}$, in contrast to existing approaches that employ polytopic methods, our approach does not provide a simple geometric representation.

\section{Finding Forward Reachable Sets Using Separating Kernels}
\label{section: main results}

Let $\mathscr{H}$ denote the Hilbert space of real-valued functions $f : \mathcal{X} \rightarrow \mathbb{R}$ on $\mathcal{X}$ equipped with the inner product $\langle \cdot, \cdot \rangle_{\mathscr{H}}$ and the induced norm $\lVert \cdot \rVert_{\mathscr{H}}$.
\begin{definition}[{\citealp{aronszajn1950theory}}]
    A Hilbert space $\mathscr{H}$ is a reproducing kernel Hilbert space if there exists a positive definite kernel function $K : \mathcal{X} \times \mathcal{X} \rightarrow \mathbb{R}$ that satisfies the following properties:
    \begin{subequations}
		\begin{align}
			& K(x, \cdot) \in \mathscr{H} && \forall x \in \mathcal{X} \\
			\label{eqn: reproducing property}
			& f(x) = \langle f, K(x, \cdot) \rangle_{\mathscr{H}} && \forall f \in \mathscr{H}, \forall x \in \mathcal{X}
		\end{align}
	\end{subequations}
	where \eqref{eqn: reproducing property} is called the reproducing property,
	and for any $x, x' \in \mathcal{X}$, we denote $K(x, \cdot)$ in the RKHS $\mathscr{H}$ as a function on $\mathcal{X}$ such that $x' \mapsto K(x, x')$.
\end{definition}
Alternatively, by the Moore-Aronszajn theorem \citep{aronszajn1950theory}, for any positive definite kernel function $K$, there exists a unique RKHS $\mathscr{H}$ with $K$ as its reproducing kernel, where $\mathscr{H}$ is the closure of the linear span of functions $\lbrace K(x, \cdot) \rbrace$.
In other words, the Moore-Aronszajn theorem allows us to define a kernel function and obtain a corresponding RKHS.

Let $K : \mathcal{X} \times \mathcal{X} \rightarrow \mathbb{R}$ be the reproducing kernel function for the RKHS $\mathscr{H}$ on $\mathcal{X}$. We stipulate that the kernel also satisfies the conditions for being a \emph{completely separating} kernel. 
This property ensures that the kernel function can be used to learn the support of any probability measure on $\mathcal{X}$.
We induce a metric $d_{K}$ on $\mathcal{X}$ via the kernel function $K$,
\begin{equation}
    \label{eqn: kernel metric}
    d_{K}(x, x') = \sqrt{K(x, x) + K(x', x') - 2K(x, x')}
\end{equation}
which is known as the kernel-induced metric on $\mathcal{X}$, and is equivalent to the Euclidean metric if $K$ is completely separating and continuous.
%
%
Following \citet{de2014learning}, we define the notion of $\mathscr{H}$ separating a subset $C \subset \mathcal{X}$.
\begin{theorem}[{\citealp[Theorem~1]{de2014learning}}]
    \label{thm: separating RKHS}
    An RKHS $\mathscr{H}$ separates a subset $C \subset \mathcal{X}$ if for all $x \not\in C$, there exists $f \in \mathscr{H}$ such that $f(x) \neq 0$ and $f(x') = 0$ for all $x' \in C$. In this case we also say that the corresponding reproducing kernel separates $C$.
\end{theorem}
According to Theorem \ref{thm: separating RKHS}, we need to determine the existence of a function in $\mathscr{H}$ that acts as a classifier for a subset $C \subset \mathcal{X}$.
Since the functions $f \in \mathscr{H}$ have the form $f = \sum_{i} \alpha_{i} K(x_{i}, \cdot)$,
this amounts to choosing a kernel function which exhibits the separating property.
Note that not all kernel functions can separate every subset $C \subset \mathcal{X}$. In order to ensure that this is possible, the notion of \emph{completely} separating kernels is defined.
\begin{definition}[{\citealp[Definition~2]{de2014learning}}]
    \label{defn: separating kernel}
    A reproducing kernel Hilbert space $\mathscr{H}$ satisfying the assumption that for all $x, x' \in \mathcal{X}$ with $x \neq x'$ we have $K(x, \cdot) \neq K(x', \cdot)$ is called completely separating if $\mathscr{H}$ separates all the subsets $C \subset \mathcal{X}$ which are closed with respect to the metric $d_{K}$. In this case, we also say that the corresponding reproducing kernel is completely separating.
\end{definition}
For example, with $\mathcal{X} = \mathbb{R}^{n}$, the Abel kernel $K(x, x') = \exp(-\lVert x - x' \rVert_{2}/\sigma)$,
where $\sigma > 0$, is a completely separating kernel function \citep[Proposition~5]{de2014learning}.
Thus, by properly selecting the kernel function $K$, we ensure that we can classify any subset $C \subset \mathcal{X}$.

We now turn to a discussion of the classifier $F$ in \eqref{eqn: forward reachable set in RKHS}.
%
%
The following theorems are reproduced with modifications for simplicity from \cite{de2014learning}, and
ensure we can define $\mathscr{F}(x_{0})$ as
the support of $x_{N}$, and that \eqref{eqn: forward reachable set in RKHS} holds.
\begin{proposition}[{\citealp[Proposition~2]{de2014learning}}]
    \label{prop: prop2}
    Assume that for all $x, x' \in \mathcal{X}$ with $x \neq x'$, $K(x, \cdot) \neq K(x', \cdot)$, the RKHS $\mathscr{H}$ with kernel $K$ is separable, and $K$ is measurable with respect to the product $\sigma$-algebra $\mathscr{B}(\mathcal{X}) \otimes \mathscr{B}(\mathcal{X})$.
    There exists a unique closed subset $\mathscr{F}(x_{0}) \subset \mathcal{X}$ with $\mathbb{P}_{N}(x_{N} \in \mathscr{F}(x_{0})) = 1$ satisfying the following property: if $C$ is a closed subset of $\mathcal{X}$ and $\mathbb{P}_{N}(x_{N} \in C) = 1$ then $\mathscr{F}(x_{0}) \subset C$.
\end{proposition}
For any subset $C \subset \mathcal{X}$, let $\mathscr{H}_{C}$ denote the closure of the linear span of functions $\lbrace K(x, \cdot) \,|\, x \in C \rbrace$, and define $P_{C} : \mathscr{H} \rightarrow \mathscr{H}$ as the orthogonal projection onto the closed subspace $\mathscr{H}_{C}$.
Define the function $F_{C} : \mathcal{X} \rightarrow \mathbb{R}$ such that $F_{C}(x) := \langle P_{C} K(x, \cdot), K(x, \cdot) \rangle_{\mathscr{H}}$.
Using this definition, we can now define the support $\mathscr{F}(x_{0})$ in terms of a function $F$.  
\begin{theorem}[{\citealp[Theorem~3]{de2014learning}}]
    Under the assumptions of Proposition \ref{prop: prop2} and the assumption that $K(x, x) = 1$ for all $x \in \mathcal{X}$,
    if $\mathscr{H}$ separates the support $\mathscr{F}(x_{0})$ of the measure $\mathbb{P}_{N}$, then $\mathscr{F}(x_{0}) = \lbrace x \in \mathcal{X} \,|\, F(x) = 1 \rbrace$.
\end{theorem}
%
Thus, we aim to learn the forward reachable set \eqref{eqn: forward reachable set in RKHS} by identifying a function $F \in \mathscr{H}$ that separates the support in $\mathcal{X}$.


\subsection{Estimating Forward Reachable Sets}

The classifier $F$ is an element of the RKHS $\mathscr{H}$, meaning it has the form $F = \sum_{i} \alpha_{i} K(x_{i}, \cdot)$ and admits a representation in terms of the finite support $\lbrace x_{i} \rbrace_{i=1}^{M} \in \mathcal{X}$, $M \in \mathbb{N}$, given by:
\begin{equation}
    \label{eqn: estimator form}
    \tilde{F}(x) = \sum_{x_{i} \in \mathcal{X}} \alpha_{i} K(x_{i}, x)
\end{equation}
where $\alpha_{i} \in \mathbb{R}$ are coefficients that depend on $x$.
Let $\mathcal{D} = \lbrace x_{N}^{i} \rbrace_{i=1}^{M}$ be a sample of terminal states at time $N$ taken i.i.d. from the system evolution.
We seek an approximation of the forward reachable set $\mathscr{F}(x_{0})$ in \eqref{eqn: forward reachable set in RKHS} using $\mathcal{D}$.
Thus, we form an estimate $\tilde{F} \in \mathscr{H}$ of $F$ using the form in \eqref{eqn: estimator form} with data $\mathcal{D}$.
%
%
We can view the estimate $\tilde{F}$ as the solution to a regularized least-squares problem.
\begin{equation}
    \label{eqn: regularized least squares}
    \min_{\tilde{F}}
    \frac{1}{M} \sum_{i=1}^{M} \vert K(x_{N}^{i}, \cdot) - \tilde{F}(x_{N}^{i}) \vert^{2} + \lambda \lVert \tilde{F} \rVert_{\mathscr{H}}^{2}
\end{equation}
where $\lambda > 0$ is the regularization parameter.
The solution to \eqref{eqn: regularized least squares} is unique and admits a closed-form solution, given by:
\begin{equation}
    \label{eqn: estimator}
    \tilde{F}(x) = \Phi^{\top} (G + M \lambda I)^{-1} \Phi
\end{equation}
where $\Phi$ is called a feature vector, with elements $\Phi_{i} = K(x_{N}^{i}, x)$, and
$G = \Phi \Phi^{\top} \in \mathbb{R}^{M \times M}$ is known as the Gram matrix, with elements $g_{ij} = K(x_{N}^{i}, x_{N}^{j})$.
A point $x \in \mathcal{X}$ is estimated to belong to the support of $x_{N}$ if $\tilde{F}(x) \geq 1 - \tau$, where $\tau$ is a threshold parameter that depends on the sample size $M$.
Thus, we form an approximation $\tilde{\mathscr{F}}(x_{0})$ of the forward reachable set $\mathscr{F}(x_{0})$ in \eqref{eqn: forward reachable set in RKHS} as:
\begin{equation}
    \label{eqn: approximate forward reachable set in RKHS}
    \tilde{\mathscr{F}}(x_{0}) = \lbrace x \in \mathcal{X} \,|\, \tilde{F}(x) \geq 1 - \tau \rbrace
\end{equation}
%
Using this representation, we obtain an estimator $\tilde{F}(x)$ which can be used to determine if a point $x \in \mathcal{X}$ is in the approximate forward reachable set $\tilde{\mathscr{F}}(x_{0})$.


\subsection{Convergence}

We characterize
the conditions for the convergence of the empirical forward reachable set $\tilde{\mathscr{F}}(x_{0})$ to $\mathscr{F}(x_{0})$
via the Hausdorff distance.
%
We assume that $\mathcal{X}$ is a topological metric space with metric $d_{K}$.
\begin{definition}[Hausdorff Distance]
    Let $A, B$ be nonempty subsets of the metric space $(\mathcal{X}, d_{K})$.
    The Hausdorff distance $d_{H}$ between sets $A$ and $B$ is defined as
    \begin{align}
        d_{H}(A, B) := \inf \lbrace \varepsilon > 0 \,|\, A \subseteq B_{\varepsilon}, B \subseteq A_{\varepsilon} \rbrace
    \end{align}
    where $A_{\varepsilon} := \bigcup_{y \in A} \lbrace x \in \mathcal{X} \,|\, d_{K}(x, y) \leq \varepsilon \rbrace$ denotes the set of all points in $\mathcal{X}$ within a ball of radius $\varepsilon$ around $A$.
\end{definition}
The Hausdorff distance gives us a method to measure convergence of the estimate $\tilde{\mathscr{F}}(x_{0})$ to the true reachable set $\mathscr{F}(x_{0})$. In fact, $\lim_{M \rightarrow \infty} d_{H}(\tilde{\mathscr{F}}(x_{0}), \mathscr{F}(x_{0})) = 0$ almost surely under mild conditions on the regularization and threshold parameters, $\lambda$ and $\tau$ \citep[Theorem~6]{de2014learning}.
As the sample size $M$ increases, if $\tau$ is chosen according to:
\begin{equation}
    \label{eqn: tau}
    \tau = 1 - \min_{1 \leq i \leq M} \tilde{F}(x_{N}^{i})
\end{equation}
and under the condition that $\lim_{M \rightarrow \infty} \lambda = 0$, the empirical forward reachable set $\tilde{\mathscr{F}}(x_{0})$ converges in probability almost surely to the true forward reachable set $\mathscr{F}(x_{0})$.

Furthermore, \citet[Theorem~7]{de2014learning} shows that the approximation admits finite sample bounds on the error of the estimated classifier function.
Let $T : \mathscr{H} \rightarrow \mathscr{H}$ be the integral operator with kernel $K$.
If $\sup_{x \in \mathcal{X}} \lVert T^{-s/2} T^{\dagger} T K(x, \cdot) \rVert \leq \infty$ with $0 < s \leq 1$, where $T^{\dagger}$ denotes the pseudo-inverse of $T$, and the eigenvalues of $T$ satisfy
$\nu_{j} = \mathcal{O}(j^{-1/b})$
for some $0 < b \leq 1$ (see \citealp[Proposition~3]{caponnetto2007optimal}), then for $M \geq 1$ and $\delta > 0$, if $\lambda = M^{-1/(2s + b + 1)}$ and $D$ is a suitable constant, then
\begin{equation}
    \sup_{x \in \mathcal{X}} \vert F(x) - \tilde{F}(x) \vert \leq D \biggl( \frac{1}{M} \biggr)^{\frac{s}{2s + b + 1}}
\end{equation}
with probability $1 - 2e^{-\delta}$.
These guarantees
ensure the accuracy of our results in a probabilistic sense, meaning that
by increasing the sample size $M$, the estimate will approach the true result, and that for a finite sample size $M$, our result is close to the true result with high probability.

As a final remark, we note that although the formulation in \eqref{eqn: forward reachable set in RKHS} is extensible to level set estimation, which involves learning a set
$\mathscr{F}(x_{0}; \alpha) := \lbrace x \in \mathcal{X} \,|\, F(x) \geq \alpha \rbrace$
with $\alpha \in [0, 1]$,
the convergence of $\tilde{F}$ to $F$ does not imply that the approximate level sets converge to the true level sets.
This requires further analysis that is beyond the scope of the current work.


\section{Numerical Results}
\label{section: numerical experiments}



For all examples, we choose an Abel kernel $K(x, x') = \exp(- \lVert x - x' \rVert_{2}/\sigma)$, $\sigma > 0$. We chose the parameters to be $\sigma = 0.1$, $\tau$ according to \eqref{eqn: tau}, and $\lambda = 1/M$, where $M$ is the sample size used to construct the classifier.
As noted in \citet{de2014learning}, $\sigma$ can be chosen via cross-validation and we require that $\lambda$ be chosen such that $\lim_{M \rightarrow \infty} \lambda = 0$.
Numerical experiments were performed in Matlab on an AWS cloud computing instance, and computation times were obtained using Matlab’s Performance Testing Framework. Code to reproduce the analysis and all figures is available at: \url{https://github.com/unm-hscl/ajthor-ortiz-L4DC2021}


\subsection{Clohessy-Wiltshire-Hill System}

\begin{figure*}
    \centering
    \includegraphics[height=150pt,keepaspectratio]{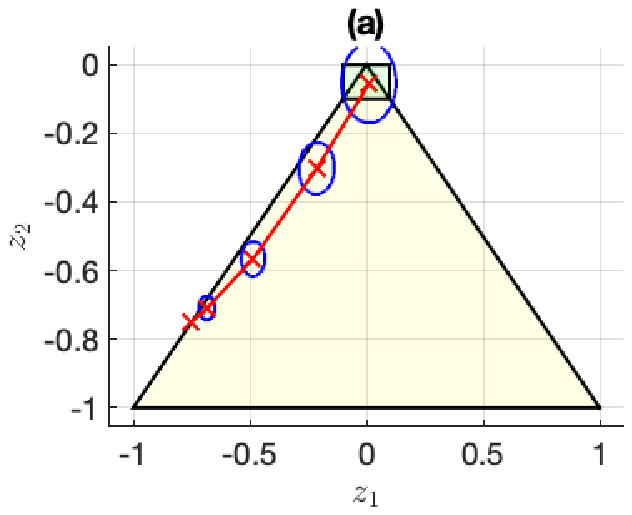}
    \includegraphics[height=150pt,keepaspectratio]{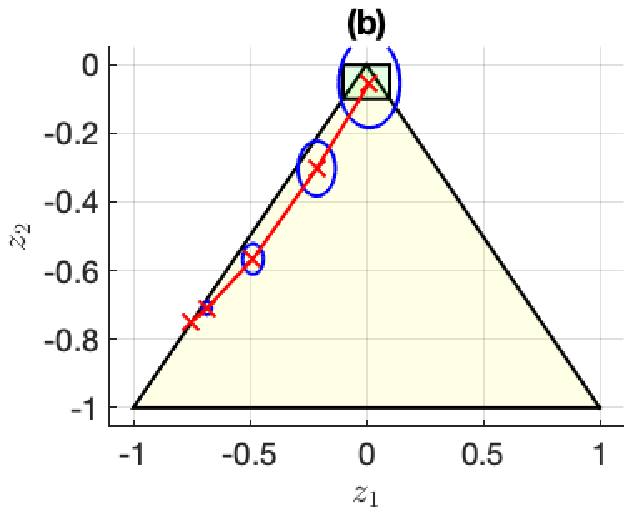}
    \caption{%
        (a) The cross-section of the approximate forward reachable set $\tilde{\mathscr{F}}(x_{0})$ computed using $M = 100$ trajectories over a time horizon of $N = 5$ from the initial condition $z_{0} = [-0.75, -0.75, 0, 0]^{\top}$ with $\dot{x} = \dot{y} = 0$ is a good approximation of
        (b) the cross-section of the forward reachable set computed using the mean and variance of the Gaussian disturbance.
        The line-of-sight cone of the spacecraft is shown in yellow, the target set is in green, and the unperturbed trajectory and initial condition is shown in red.%
    }
    \label{fig: cwh reachable sets}
\end{figure*}

We first consider a realistic example of
spacecraft rendezvous and docking to compare our technique against a known result.
The dynamics of a CWH system are defined in \citet{lesser2013stochastic} as:
\begin{align}
    \label{eqn: cwh dynamics}
    \begin{split}
        \begin{aligned}[c]
        	\ddot{x} - 3\omega^{2} x - 2 \omega \dot{y} &= F_{x}/m_{d} &&&
        	\ddot{y} + 2 \omega \dot{x} &= F_{y}/m_{d}
        \end{aligned}
    \end{split}
\end{align}
with state $z = [x, y, \dot{x}, \dot{y}]^{\top} \in \mathcal{X} \subseteq
\mathbb{R}^{4}$, input $u = [F_{x}, F_{y}]^{\top} \in \mathcal{U} \subseteq \mathbb{R}^{2}$, where
$\mathcal{U} = [-0.1, 0.1] \times [-0.1, 0.1]$, and parameters $\omega$, $m_d$.
We first discretize the dynamics in time and then apply an affine Gaussian
disturbance $w$, which is a stochastic process defined on $(\mathcal{W}, \mathscr{B}(\mathcal{W}))$ with variance $\Sigma =
\textnormal{diag}(1 \times 10^{-4}, 1 \times 10^{-4}, 5 \times 10^{-8}, 5 \times
10^{-8})$ such that $w_{k} \sim \mathcal{N}(0, \Sigma)$.

With initial condition $z_{0} = [-0.75, -0.75, 0, 0]^{\top}$, we compute an open-loop control policy $\pi$ using a chance-constrained algorithm from \citet{vinod2019affine} implemented in SReachTools \citep{SReachTools}.
The control policy is designed to dock with another spacecraft while remaining within a line of sight cone.
We then simulated $M = 100$ trajectories
to collect a sample $\mathcal{D}$ consisting of the resulting terminal states $\lbrace z_{N}^{i} \rbrace_{i=1}^{M}$.
A classifier was then computed using \eqref{eqn: estimator}.
In order to depict the sets graphically, we then chose a small region around the mean of the observed state values and computed a visual representation of the forward reachable set by sampling 10,000 points uniformly over the region and connecting evaluations where $\tilde{F}(\cdot) = 1$ along the set boundary.


Figure \ref{fig: cwh reachable sets} compares the computed approximate forward reachable sets
to forward reachable sets computed
using the known Gaussian disturbance properties. 
We can see that the approximate forward reachable sets $\tilde{\mathscr{F}}(x_{0})$ in Figure \ref{fig: cwh reachable sets}(a) are close to the variance ellipses produced from the known disturbance properties in \ref{fig: cwh reachable sets}(b), demonstrating that
our approach provides a good estimate of the support of the underlying distribution.
Computation time for the approximate forward reachable sets was 1.582 seconds using the kernel based approach.
%
%
%
This technique could be accelerated with approximative speedup techniques that are common to kernel methods, such as random Fourier features \citep{rahimi2008random} or FastFood \citep{le2013fastfood}.
Additionally, incorporating active or adversarial sampling, as in \citet{lew2020sampling}, could reduce computation time by reducing the number of observations needed to compute the classifier.


\subsection{Neural Network Verification}

We now demonstrate our approach on a set of dynamical system benchmarks with neural network controllers as described in \citet{dutta2019sherlock}.
We define a feed-forward neural network controller $\pi : \mathcal{X} \rightarrow \mathcal{U}$ as
$\pi(z) = g_{L} \circ \cdots \circ g_{0}(z, \theta_{0})$,
with activation functions $g_{i}(z, \theta_{i})$, $i = 1, \ldots, L - 1$ that depend on the parameters $\theta_{i}$.
The function $g_{0}$ maps the states into the first layer, and $g_{L}$ is a function that maps the output of the last layer into the control space.
%
After generating
observations of the states, we then assume no knowledge of the structure of $\pi$ or the dynamics.


\subsubsection{TORA Model}

\begin{figure*}
    \centering
    \includegraphics[height=150pt,keepaspectratio]{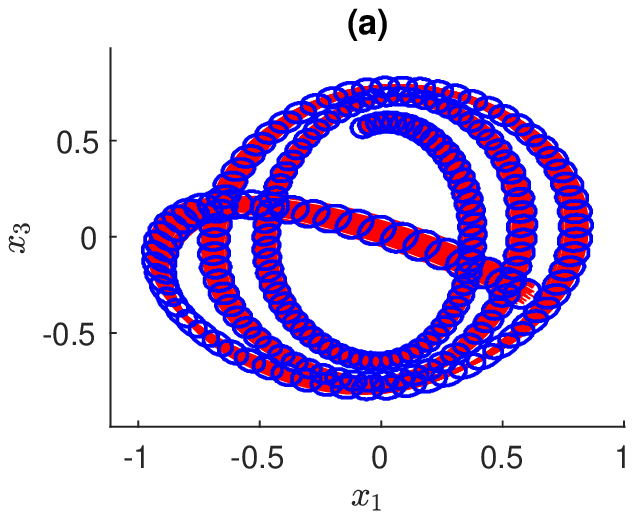}
    \includegraphics[height=150pt,keepaspectratio]{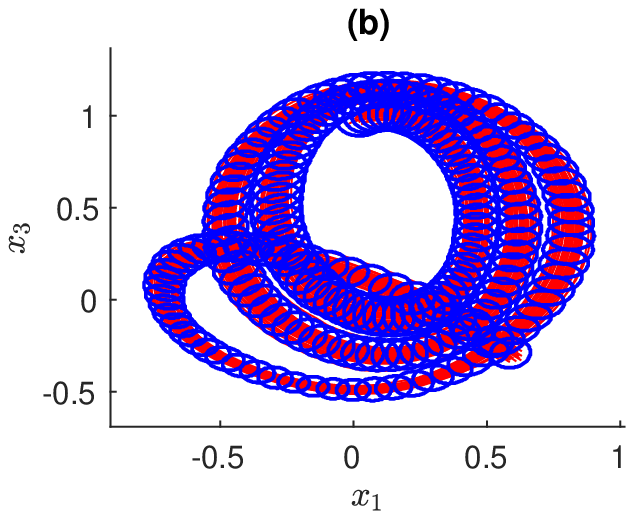}
    \caption{%
        (a) Cross-section of the approximate forward reachable set $\tilde{\mathscr{F}}(x_{0})$ for the TORA model over a time horizon of $N = 200$ using $M = 50$ trajectories, validated against simulated trajectories.
        (b) Cross-section of the approximate forward reachable sets $\tilde{\mathscr{F}}(x_{0})$ for the TORA model with an affine beta distribution disturbance $w_{k} \sim 0.01 \; \mathrm{Beta}(2, 0.5)$.
        The observed trajectories are shown in red.%
    }
    \label{fig: tora}
\end{figure*}

Consider a translational oscillations by a rotational actuator (TORA) model from \citet{jankovic1996tora} with a neural network controller \citep[Benchmark~9]{dutta2019sherlock}, with dynamics given by:
\begin{align}
    \label{eqn: tora dynamics}
    \begin{split}
        \begin{aligned}[c]
            \dot{x}_{1} &= x_{2} &&&
            \dot{x}_{3} &= x_{4} \\
            \dot{x}_{2} &= -x_{1} + 0.1 \sin(x_{3}) &&&
            \dot{x}_{4} &= u
        \end{aligned}
    \end{split}
\end{align}
where $u$ is the control input, chosen by the neural network controller. The neural network is trained via an MPC algorithm to keep all state variables within the range $[-2, 2]$.
We presume an initial distribution that is uniform over the range $[0.6, 0.7] \times [-0.7, -0.6] \times [-0.4, -0.3] \times [0.5, 0.6]$ and collected a sample from $M = 50$ simulated trajectories over a time horizon of $N = 200$.
Using \eqref{eqn: estimator}, we then computed a classifier to indicate whether a given point is within the approximate forward reachable set $\tilde{\mathscr{F}}(x_{0})$.
As before, we chose a small region around the mean of the observed trajectories and
created a visual representation of
the approximate forward reachable set
by connecting evaluations where $\tilde{F}(\cdot)=1$ along the set boundary.
We validate our approach via Monte Carlo simulation.
Figure \ref{fig: tora}(a) shows a cross-section of the approximate forward reachable set $\tilde{\mathscr{F}}(x_{0})$ \eqref{eqn: approximate forward reachable set in RKHS}.
As expected, we can see the forward reachable set encompasses the simulated trajectories.  
The computation time was 52.907 seconds for computing the approximate forward reachable sets $\tilde{\mathscr{F}}(x_{0})$ over the time horizon $N = 200$.


We then add an arbitrarily chosen affine disturbance $w$ to the dynamics in \eqref{eqn: tora dynamics}, which is a stochastic process consisting of the random variables
with a beta distribution
$w_{k} \sim 0.01 \; \mathrm{Beta}(\alpha, \beta)$, with PDF
$f(x \,|\, \alpha, \beta) =
x^{\alpha-1} (1-x)^{\beta-1}/B(\alpha, \beta)$
where $B(\alpha, \beta) = \Gamma(\alpha) \Gamma(\beta)/\Gamma(\alpha + \beta)$ and $\Gamma$ is the Gamma function, with shape parameters $\alpha = 2$, $\beta = 0.5$.
We then implemented our method on the stochastic system and computed the approximate reachable sets. As shown in Figure \ref{fig: tora}(b), the cross-section shows a larger variation in the trajectories due to the added disturbance, resulting in larger approximate forward reachable sets.




\subsubsection{Drone Model}

\begin{figure*}
    \centering
    \includegraphics[height=150pt,keepaspectratio]{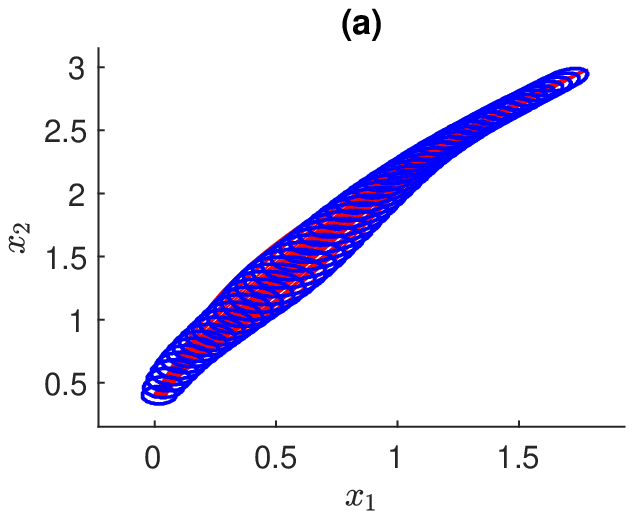}
    \includegraphics[height=150pt,keepaspectratio]{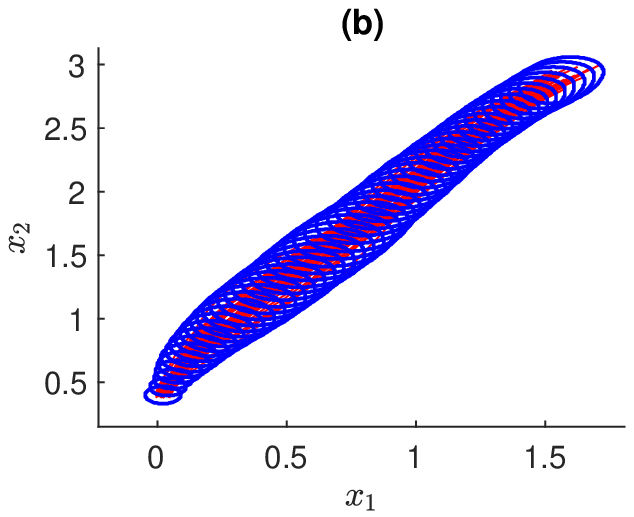}
    \caption{%
        (a) Cross-section of the approximate forward reachable sets $\tilde{\mathscr{F}}(x_{0})$ for the drone model over a time horizon $N = 50$ using $M = 50$ trajectories.
        (b) Cross-section of the approximate forward reachable sets $\tilde{\mathscr{F}}(x_{0})$ for the drone model with an affine Gaussian disturbance $w_{k} \sim \mathcal{N}(0, \Sigma)$, $\Sigma = 0.0025 I$.
        The observed trajectories are shown in red.%
    }
    \label{fig: drone}
\end{figure*}

Lastly, consider a 12-DOF quadrotor system with nonlinear dynamics defined in \citet{bansal2016learning} with a neural network controller described in \citet{dutta2019sherlock}.
This system has proven challenging for existing reachability tools, since the trajectories diverge locally before converging, meaning the forward reachable set is difficult to compute using over-approximative interval based methods.
%
%
The system is controlled via four inputs
$u \in \mathbb{R}^{4}$ which are determined by a neural network controller.
Following \citet{dutta2019sherlock}, we chose an arbitrary state close to the origin $z \in \mathbb{R}^{12}$
with initial distribution uniform over the range
$x_{0} \sim z + 0.05 \; U(0, 1)$
and collected a sample from $M = 50$ simulated trajectories over a time horizon of $N = 50$.
We calculated the approximate forward reachable sets $\tilde{\mathscr{F}}(x_{0})$ as before.
Figure \ref{fig: drone}(a) shows a cross-section of the approximate forward reachable sets for the problem.
As expected, we can see the approximate forward reachable sets $\tilde{\mathscr{F}}(x_{0})$ \eqref{eqn: approximate forward reachable set in RKHS} are centered around the trajectories (in red), indicating that the approximate forward reachable sets are a good approximation of the support.
The computation time for calculating the forward reachable sets was 12.830 seconds using the kernel method based approach.

We further demonstrate the capability of our algorithm by applying an arbitrarily chosen affine disturbance $w$ to the dynamics, which as before is a stochastic process consisting of the random variables $w_{k} \sim \mathcal{N}(0, \Sigma)$, with $\Sigma = 0.0025 I$.
Figure \ref{fig: drone}(b) depicts the approximate forward reachable sets $\tilde{\mathscr{F}}(x_{0})$ for the stochastic nonlinear system.
This shows that our method can handle high-dimensional, stochastic nonlinear systems controlled by neural networks, and compute the approximate forward reachable sets with efficient computational time.


\section{Conclusion \& Future Work}
\label{section: conclusion}

We presented a method for computing forward reachable sets using reproducing kernel Hilbert spaces with separating kernel functions. We demonstrated the effectiveness of the technique at performing neural network verification, and validated its accuracy on a satellite rendezvous and docking problem with Clohessy-Wiltshire-Hill dynamics.
This technique is scalable, computationally efficient, and model-free, making it well-suited for systems with learning-enabled components.
We plan to investigate extension of this to other reachability problems, such as safety problems that require calculation of the backward reachable set. 


\acks{This material is based upon work supported by the National Science Foundation under NSF Grant Number CNS-1836900.  Any opinions, findings, and conclusions or recommendations expressed in this material are those of the authors and do not necessarily reflect the views of the National Science Foundation.  This research was supported in part by the Laboratory Directed Research and Development program at Sandia National Laboratories, a multimission laboratory managed and operated by National Technology and Engineering Solutions of Sandia, LLC., a wholly owned subsidiary of Honeywell International, Inc., for the U.S. Department of Energy’s National Nuclear Security Administration under contract DE-NA-0003525.  The views expressed in this article do not necessarily represent the views of the U.S. Department of Energy or the United States Government.}


\bibliography{bibliography}

\end{document}